\documentclass[12pt]{article}
\usepackage{amsmath}
\usepackage{amsfonts}
\usepackage{amscd}
\usepackage[english]{babel}
\usepackage[cp866]{inputenc}

\newtheorem{lemma}{Lemma}
\newtheorem{corollary}{Corollary}
\newtheorem{proposition}{Proposition}

\newtheorem{definition}{Definition}
\newtheorem{remark}{Remark}

\def\g {\mathfrak{g}}
\def\h {\mathfrak{h}}
\def\p {\mathfrak{p}}
\title{$U(n+1)\times U(p+1)$ - invariant Hermitian metrics
with Hermitian tensor Ricci on the manifold $S^{2n+1}\times
S^{2p+1}$.}
\author{Daurtseva N. A.,\ Kemerovo State
University}

\begin{document}
\maketitle
\begin{abstract}
Invariant complex structures on the homogeneous manifold
$U(n+1)/U(n)\times U(p+1)/U(p)$ are reseached. The critical point
of the functional of the scalar curvature is found.
\end{abstract}
Let consider the product $S^{2n+1}\times S^{2p+1}$
as homogeneous space $U(n+1)/U(n)\times U(p+1)/U(p)$. Suppose,
that $n$ and $p$ are not equal to zero simultaneously. Denote
$\g_1$ and $\h_1$ ($\g_2$ and $\h_2$) Lie algebras of Lie groups
$U(n+1)$ and $U(n)$ ($U(p+1)$ and $U(p)$) respectively. As group
$U(n)$ is embedded into $U(n+1)$ by usual way, then $\h_j$ is
embedded into $\g_j$ by the following way:
$$
M\in\h_j\mapsto \left(\begin{array}{cc} 0 & 0\\
0 & M\end{array}\right)\in\g_j,
$$
where $j=1,2$. Let define the basis in $\g_1\times\g_2$. Let
$E^j_{\nu\mu}$ is matrix with 1 on the $(\nu,\mu)$-place and
other zero elements. Define:
$$
Z^j_{\nu\mu}=E^j_{\nu\mu}-E^j_{\mu\nu},\
T^j_{\nu\mu}=E^j_{\nu\mu}+E^j_{\mu\nu},\ 0\leq\mu <\nu\leq n,\
j=1,2
$$
Matrix $Z^j_{\nu\mu},\ iT^j_{\nu\mu}$ (where $j=1,2$) form basis
of product $\g_1\times\g_2$. Take decomposition $\g_j=\h_j\oplus
\p_j$, where $\mathfrak{p}_j$ has basis $X^j=\frac12iT^j_{00}$,
$Y^j_{2\nu-1}=Z^j_{\nu 0}$, $Y^j_{2\nu}=iT^j_{\nu 0}$. So,
manifold $S^{2n+1}\times S^{2p+1}$ viewed as homogeneous space
$U(n+1)/U(n)\times U(p+1)/U(p)$ has basis
$X^1,Y^1_{2\nu-1},Y^1_{2\nu},X^2,Y^2_{2\mu-1},Y^2_{2\mu}$,
$1\leq\nu\leq n$, $1\leq\mu\leq p$.
\begin{proposition}
$$
\begin{array}{ll}
1)\ [\mathfrak{p}_0,\mathfrak{p}_1]\subset\mathfrak{p}_1: &
[X^1,Y^1_{2\nu-1}]=-Y^1_{2\nu}, [X^1,Y_{2\nu}^1]=Y_{2\nu-1}^1,\\
2)\
[\mathfrak{p}_1,\mathfrak{p}_1]\subset\h_1\oplus\mathfrak{p}_0:
& [Y^1_{2\nu-1},Y^1_{2\nu}]=-2X^1+iT^1_{\nu\nu},\\
 & [Y_{2\nu}^1,Y_{2\mu}^1]=-Z^1_{\nu\mu},\\
 & [Y_{2\nu-1}^1,Y_{2\mu-1}^1]=-Z^1_{\nu\mu},\\
 & [Y_{2\nu}^1,Y_{2\mu-1}^1]=-T^1_{\nu\mu},\\
3)\ [\mathfrak{p}_2,\mathfrak{p}_3]\subset\mathfrak{p}_3: &
[X^2,Y^2_{2\nu-1}]=-Y^2_{2\nu}, [X^2,Y_{2\nu}^2]=Y_{2\nu-1}^2,\\
4)\
[\mathfrak{p}_3,\mathfrak{p}_3]\subset\h_2\oplus\mathfrak{p}_2:
& [Y^2_{2\nu-1},Y^2_{2\nu}]=-2X^2+iT^2_{\nu\nu},\\
 & [Y_{2\nu}^2,Y_{2\mu}^2]=-Z^2_{\nu\mu},\\
 & [Y_{2\nu-1}^2,Y_{2\mu-1}^2]=-Z^2_{\nu\mu},\\
 & [Y_{2\nu}^2,Y_{2\mu-1}^2]=-T^2_{\nu\mu},\\
5)\
[\mathfrak{p}_0\oplus\mathfrak{p}_1,\mathfrak{p}_2\oplus\mathfrak{p}_3]=0
\end{array}
$$
\end{proposition}
{\bf Proof.} Proposition follows from definition of vectors
$X^i,Y^i_{2\nu-1},Y^i_{2\nu}$ ($i=1,2$).
\begin{definition}
Almost complex structure on the manifold $M$ is smooth field of
endomorphisms $J_x:T_xM\longrightarrow T_xM$, such that
$J^2_x=-Id_x$, $\forall x\in M$, where $Id_x$ is identical
endomorphism $T_xM$.
\end{definition}
Recall some known construction of complex structure on
$S^{2n+1}\times S^{2p+1}\quad$ \cite{4}. It is known that
$S^{2n+1}\times S^{2p+1}$ is principal  $S^1\times S^1$ bundle
over $\mathbb{CP}^n\times\mathbb{CP}^p$. The space
$\mathbb{CP}^n\times\mathbb{CP}^p$ and fiber $S^1\times S^1$ are
complex manifolds. If we fix complex structures on the base and
fiber, then we can choose holomorphic transition functions to get
complex structure on $S^{2n+1}\times S^{2p+1}$. All those
structures form two parametric family $I(a,c)$ ($c>0$), they are
$U(n+1)\times U(p+1)$ - invariant \cite{4}.

Consider projection:
$$
U(n+1)/U(n)\times U(p+1)/U(p)\longrightarrow U(n+1)/(U(n)\times
U(1))\times U(p+1)/(U(p)\times U(1))
$$
Obviously, that $U(n+1)/(U(n)\times U(1))\times U(p+1)/(U(p)\times
U(1))$ is product of complex projective spaces
$\mathbb{CP}^n\times\mathbb{CP}^p$, and vectors $X^1,\ X^2$ are
tangent to fiber. $I(a,c)$ acts on these vectors as
$$
I(a,c)X^1=\frac{a}{c}X^1+\frac1cX^2,\quad
I(a,c)X^2=-\frac{a^2+c^2}{c}X^1-\frac{a}{c}X^2
$$
$$
I(a,c)Y^1_{2\nu-1}=Y^1_{2\nu},\quad I(a,c)Y^2_{2\mu-1}=Y^2_{2\mu},
$$
where parameters $a$ and $c$ are real, $c>0$. As $I(a,c)$ are
$U(n+1)\times U(p+1)$ - invariant, then they defined by action of
$I(a,c)$ on the basis of the space $\p_1\times\p_2$. Denote
$\p_1\times\p_2$ as $\p$, and $\h_1\times\h_2$ as $\h$.
\begin{definition}
Almost complex structure $J$ on the manifold $M$ is called
positive associated with  2-form $\omega$, if:\\
1) $\omega(JX,JY)=\omega(X,Y)$, for all $X,Y\in TM$\\
2) $\omega(X,JX)>0$, for all nonzero $X\in TM$
\end{definition}
Fix non-degenerate invariant 2-form $\omega$:
$$
\omega=X^1\wedge X^2+\sum_{\nu = 1}^nY^1_{2\nu-1}\wedge
Y^1_{2\nu}+\sum_{\nu = 1}^pY^2_{2\nu-1}\wedge Y^2_{2\nu}
$$
on $S^{2n+1}\times S^{2p+1}$
\begin{lemma}
All complex structures $I(a,c)$ are positive associated with
$\omega$.
\end{lemma}
{\bf Proof.} For $I(a,c)$ properties 1) and 2) of definition 2 are
obvious.
\begin{corollary}
Every complex structure $I(a,c)$ defines unique $\omega$ -
associated metric
$$
g(a,c)(X,Y)=\omega(X,I(a,c)Y)
$$
\end{corollary}
These associated metrics are:
$$
g(a,c)(X^1,X^1)=1/c,\ g(a,c)(X^2,X^2)=(a^2+c^2)/c,\
g(a,c)(X^1,X^2)=-a/c
$$
$$
g(a,c)(Y^1_j,Y^1_j)=g(a,c)(Y^2_k,Y^2_k)=1,\ 1\leq j\leq 2n,\ 1\leq
k\leq 2p
$$
$$
g(a,c)(X,Y)=0,\mbox{ for other basis vectors }X\mbox{ and }Y
$$
Each metric of this family is $I(a,c)$ - Hermitian, so we obtain
two-parametric family of Hermitian manifolds $(S^{2n+1}\times
S^{2p+1},$ $g(a,c), I(a,c), \omega)$. Invariant metric induces
scalar product on $\p$.
\begin{proposition}
Invariant Riemmanian connection for $g(a,c)$ on the
$S^{2n+1}\times S^{2p+1}$ is given by formula
$D_XY=\frac12[X,Y]_{\p}+U(X,Y)$, where $U$ is symmetric bilinear
mapping $U:\p\times\p\longrightarrow\p$:
$$
U(X^1,Y^1_{2\nu-1})=\frac{2-c}{2c}\ Y^1_{2\nu},\
U(X^1,Y^1_{2\nu})=-\frac{2-c}{2c}\ Y^1_{2\nu-1},
$$
$$
U(X^1,Y^2_{2\nu-1})=-\frac{a}{c}\ Y^2_{2\nu},\
U(X^1,Y^2_{2\nu})=\frac{a}{c}\ Y^2_{2\nu-1},
$$
$$
U(X^2,Y^1_{2\nu-1})=-\frac{a}{c}\ Y^1_{2\nu},\
U(X^2,Y^1_{2\nu})=\frac{a}{c}\ Y^1_{2\nu-1},
$$
$$
U(X^2,Y^2_{2\nu-1})=\left(\frac{a^2+c^2}{c}-\frac12\right)\
Y^2_{2\nu},\
U(X^2,Y^2_{2\nu})=-\left(\frac{a^2+c^2}{c}-\frac12\right)\
Y^2_{2\nu-1},
$$
$U(X,Y)=0$ for other basis vectors $X$ and $Y$.
\end{proposition}
{\bf Proof.} Find $U$ by formula:
$2g(U(X,Y),Z)=g([Z,X]_{\p},Y)+g(X,[Z,Y]_{\p})$

\begin{proposition}
Two-parametric family of metrics $g(a,c)$ has following
characteristics:\\
1. Ricci curvature:
$$
Ric(a,c)(X^1,X^1)=2\frac{n+pa^2}{c^2},\
Ric(a,c)(X^2,X^2)=2\frac{na^2+p(a^2+c^2)^2}{c^2},
$$
$$
Ric(a,c)(X^1,X^2)=-2\frac{a}{c^2}(n+p(a^2+c^2)),
$$
$$
Ric(a,c)(Y^1_j,Y^1_j)=2(1+n-\frac1c),\ 1\leq j\leq 2n,
$$
$$
Ric(a,c)(Y^2_k,Y^2_k)=2(1+p-\frac{a^2+c^2}{c}),\ 1\leq k\leq 2p,
$$
$$
Ric(a,c)(X,Y)=0,\mbox{ for other basis vectors }X\mbox{ and }Y
$$
Proper values of Ricci curvature $\tilde{r}_i$ are
$\tilde{r}_{1,2}=\frac{x+y\pm\sqrt{(x-y^2+4z^2)}}{2}$, where
$x=2\frac{n+pa^2}{c^2}$, $y=2\frac{na^2+p(a^2+c^2)^2}{c^2}$,
$z=-2\frac{a}{c^2}(n+p(a^2+c^2))$;
$\tilde{r}_3=\tilde{r}_4=\dots=\tilde{r}_{2n+2}=2(1+n-\frac1c)$,
$\tilde{r}_{2n+3}=\tilde{r}_{2n+4}=\dots=\tilde{r}_{2n+2p+2}=2(1+p-\frac{a^2+c^2}{c})$.\\
2. Scalar curvature:
$$
s=4n\left(1+n-\frac{1}{2c}\right)+4p\left(1+p-\frac{a^2+c^2}{2c}\right)
$$
\end{proposition}
{\bf Proof.}1. Calculate Ricci curvature by formula \cite{1}:
$$Ric(X,X)=-\frac12\sum_i|[X,Z_i]_{\p}|^2-\frac12\sum_ig([X,[X,Z_i]_{\p}]_{\p},Z_i)$$
$$-\sum_ig([X,[X,Z_i]_{\h}]_{\p},Z_i)+\frac14\sum_{i,j}g([Z_i,Z_j]_{\p},X)^2-g([Z,X]_{\p},X),$$
where $Z=\sum_iU(Z_i,Z_i)$ and $Z_i$ is orthonormal basis of the
space $(\p,g)$.

In our case, orthonormal basis with respect to $g(a,c)$ is:
$Z_0=\sqrt{c}$, $Z_{i}=Y^1_{i}$, for $i=1,\dots 2n$,
$Z_{2n+1}=\frac{a}{\sqrt{c}}X^1+\frac{1}{\sqrt{c}}X^2$,
$Z_{2n+1+i}=Y^2_i$, where $i=1,\dots,2p$. Obviously, that $Z=0$.
2. Scalar curvature is calculated as trace of Ricci tensor:
$s=Ric_{ij}g^{ij}$, where $g^{ij}$ are components of $g(a,c)^{-1}$
($i,j=1,\dots,2n+2p+2$).

The family of complex structures $I(a,c)$ on $S^{2n+1}\times
S^{2p+1}$ consists of all $U(n+1)\times U(p+1)$ - invariant
almost complex structures. So, if $\mathcal{A}^+_{\omega}$ is
space of invariant almost complex structures, which are positive
associated with $\omega$, and $\mathcal{AM}_{\omega}^+$ is the
space of positive associated metrics, then:
$$
\mathcal{A}^+_{\omega}=\{I(a,c):c>0\}\qquad\mathcal{AM}^+_{\omega}=\{g(a,c):c>0\}
$$
The functional of scalar curvature is defined on the
$\mathcal{AM}^+_{\omega}$:
$$
s:\mathcal{AM}^+_{\omega}\longrightarrow\mathbb{R},\qquad
s(g)=4n(1+n-\frac{1}{2c})+4p(1+p-\frac{a^2+c^2}{2c})
$$
It is known  (see, for example \cite{2}), that critical points of
this functional on $\mathcal{AM}^+_{\omega}$ give metrics with
$I$- Hermitian Ricci tensor.
\begin{proposition}
If $n$ or $p$ is equal to zero, then there are not $U(n+1)\times
U(p+1)$ - invariant metrics $g(a,c)$ with  Hermitian Ricci tensor
on $S^{2n+1}\times S^{2p+1}$. If $n$ and $p$ are not equal to
zero, then metric $g(a,c)$, when $a=0$, $c=\sqrt{\frac{n}{p}}$ has
$I(a,c)$ - Hermitian Ricci tensor.
\end{proposition}
{\bf Proof.} Find partial derivatives  of $s(a,c)$ with respect to
$a$ and $c$:
$$
\frac{\partial s}{\partial a}=-p\frac{a}{c}
$$
$$
\frac{\partial s}{\partial c}=\frac{n-p(c^2-a^2)}{2c^2}
$$
So, if $n$ or $p$ is equal to zero, then $s$ has no critical
points. If $n$ and $p$ are not equal to zero, then functional $s$
takes maximal value $4n(n+1)+4p(1+p)-4\sqrt{np}$ at point $a=0$,
$c=\sqrt{\frac{n}{p}}$.

\begin{remark}
One can shows, that if $n=p$, $a=0$, $c=\sqrt{\frac{n}{p}}=1$,
then:
$$
Ric=2ng
$$
Therefore the above metric is Einstein. If $n\neq p$, then above
metrics are not  Einstein.
\end{remark}
Let $n\leq p$, we can apply result of \cite{3}:
\begin{proposition}
Sectional curvature of metric $g(0,\sqrt{\frac{n}{p}})$ satisfied
to the following inequalities: \\
1. If $0<\frac{n}{p}\leq\frac19$, then $4-3\sqrt{\frac{p}{n}}\leq
K\leq\sqrt{\frac{p}{n}}$. Minimal value is obtained on the
bivector $Y^1_{2l-1}\wedge Y^1_{2l}$ ($l=1,\dots,n$), and maximal
on the $\sqrt{c}X^1\wedge Y^1_{i}$ ($i=1,\dots,2n$).\\
2. If $\frac19<\frac{n}{p}\leq\frac{9}{16}$, then
$4-3\sqrt{\frac{p}{n}}\leq K\leq 4-3\sqrt{\frac{n}{p}}$. Minimal
value is obtained on the bivector $Y^1_{2l-1}\wedge Y^1_{2l}$
($l=1,\dots,n$), and maximal on the $Y^2_{2m-1}\wedge Y^2_{2m}$ ($m=1,\dots,p$).\\
3. If $\frac{9}{16}<\frac{n}{p}\leq 1$, then $0\leq K\leq
4-3\sqrt{\frac{n}{p}}$. Minimal value is obtained on bivectors
$X^1\wedge X^2$, $Y^1_{2l-1}\wedge Y^2_{2m-1}$ and $Y^1_{2l}\wedge
Y^2_{2m}$ ($l=1,\dots,n$, $m=1,\dots,p$), maximal on the
$Y^2_{2l-1}\wedge Y^2_{2l}$.
\end{proposition}

\end{document}